\documentclass[11pt]{article}
\newtheorem{theorem}{Theorem}
\newtheorem{lemma}[theorem]{Lemma}

 \newcommand{\nl}{\newline}
 \newcommand{\dist}{{\rm dist}}
 \newcommand{\N}{{\bf N}}
 
\newcommand{\R}{{\bf R}}

 \newcommand{\diver}{{\rm div}}

 \newcommand{\codim}{{\rm codim}}

\newcommand{\ia}{({\rm i})}
\newcommand{\ib}{({\rm ii})}
\newcommand{\ic}{({\rm iii})}

\title{On a class of Rellich inequalities}

\author{
G. Barbatis\footnote{Department of Mathematics,
 University of Ioannina, 45110 Ioannina, Greece}
 \and A. Tertikas
 \footnote{Department of Mathematics,
 University of Crete, 71409 Heraklion, Greece and  \nl
Institute of  Applied and Computational Mathematics,
FORTH, 71110 Heraklion, Greece}
}

\begin{document}

\maketitle

\begin{center}
{\em Dedicated to Professor E.B. Davies on the occasion of his 60th birthday}
\end{center}

\

\begin{abstract}

We prove Rellich and improved Rellich inequalities that involve
the distance function from a hypersurface of codimension $k$,
under a certain geometric assumption. In case the distance is taken
from the boundary, that assumption is the convexity of the domain.
We also discuss the best constant of these inequalities.

\noindent {\bf AMS Subject Classification: 35J20 (35P20, 35P99,
26D10, 47A75)}\nl {\bf Keywords: Rellich inequality, best
constants, distance function }
\end{abstract}

\section{Introduction}

The classical Rellich inequality states that for $p>1$
\begin{equation}
\int_0^{\infty}|u''|^pdt \geq\frac{(p-1)^p(2p-1)^p}{p^{2p}}\int_0^{\infty}\frac{|u|^p}{t^{2p}}dt\, ,
\label{man}
\end{equation}
for all $u\in C^{\infty}_c(0,\infty)$. A multi-dimensional version of (\ref{man}) for $p=2$
is also classical and states that for any $\Omega\subset\R^N$, $N\geq 5$, there holds
\begin{equation}
\int_{\Omega}(\Delta u)^2dx \geq\frac{N^2(N-4)^2}{16}\int_{\Omega}\frac{u^2}{|x|^{4}}dx\, ,
\label{man1}
\end{equation}
for all $u\in C^{\infty}_c(\Omega)$.

Davies and Hinz \cite{DH} generalized (\ref{man1}) and showed that for any $p\in (1,N/2)$ there holds
\begin{equation}
\int_{\Omega}|\Delta u|^pdx \geq\biggl(\frac{(p-1)N|N-2p|}{p^2}\biggr)^p
\int_{\Omega}\frac{|u|^p}{|x|^{2p}}dx\; , \qquad
u\in C^{\infty}_c(\Omega\setminus\{0\})\, .
\label{r1}
\end{equation}
Inequality (\ref{man}) has also been generalized to higher dimensions in another
direction, where the singularity involves the distance $d(x)=\dist(x,\partial\Omega)$.
Owen \cite{O} proved among other results that if $\Omega$ is bounded and convex then
\begin{equation}
\int_{\Omega}(\Delta u)^2dx \geq \frac{9}{16}\int_{\Omega}\frac{u^2}{d(x)^{4}}dx\; , \qquad
u\in C^{\infty}_c(\Omega)\, .
\label{r2}
\end{equation}
Recently,   an improved version of (\ref{man1}) has been
established in \cite{TZ}. Among several other results they showed
that for a bounded domain $\Omega$ in $\R^N$, $N \geq 5$, there holds
\begin{equation}
\int_{\Omega}(\Delta u)^2dx \geq
\frac{N^2(N-4)^2}{16}\int_{\Omega}\frac{u^2}{|x|^{4}}dx
+\biggl(1+\frac{N(N-4)}{8}\biggr)\sum_{i=1}^{\infty}\int_{\Omega}
\frac{u^2}{|x|^4}X_1^2X_2^2\ldots X_i^2 dx\, , \label{xxa}
\end{equation}
as well as
\begin{eqnarray*}
\int_{\Omega}(\Delta u)^2dx\geq
\frac{N^2}{4}\int_{\Omega}\frac{|\nabla u|^2}{|x|^2}dx
+\frac{1}{4}\sum_{i=1}^{\infty}\int_{\Omega}
\frac{|\nabla u|^2}{|x|^2}X_1^2X_2^2\ldots X_i^2 dx \; , \\
\end{eqnarray*}
 for all $u\in C^{\infty}_c(\Omega\setminus\{0\})$. Here
$X_k$ are iterated logarithmic functions; see (\ref{va}) for
precise the definition.

Rellich inequalities have various applications in the study of
fourth-order elliptic and parabolic PDE's; see e.g. \cite{DH,O,B}.
Improved Rellich inequalities are useful if critical potentials
are additionally present. As a simplest example, one obtains
information on the existence of solution and asymptotic behavior for
the equation $u_t=-\Delta^2+V$ for critical potentials $V$. 
Corresponding problems for improved Hardy's inequalities have recently attracted considerable
attention: see \cite{BV,BM,BFT1} and references therein.

Our aim in this paper is to obtain sharp improved versions of inequalities (\ref{r1}) and (\ref{r2}),
where additional non-negative terms are present in the respective right-hand sides.
At the same time we obtain some new improved Rellich inequalities which are new even at the
level of plain Rellich inequalities; these involve the distance to a surface $K$ of intermediate codimension.

\

{\bf\large Statement of results}

Before stating our theorems let us first introduce some notation.
We denote by $\Omega$ a domain in $\R^N$, $N \geq 2$. For the sake of simplicity all functions considered below
are assumed to be real-valued; in relation to this we note however that minor modifications of the proofs or a suitable application
of \cite[Lemma 7.5]{D} can yield the validity of Theorems 1-3 below for complex-valued functions $u$.
We let $K$ be a closed,
piecewise smooth surface of codimension $k$, $k=1,\ldots,N$. We do not assume that
$K$ is connected but only that it has finitely many connected components. In the case $k=N$
we assume that $K$ is a finite union of points while in the case $k=1$ we assume that $K=\partial\Omega$.
We then set
\[
d(x) = \dist(x, K) \, ,
\]
and assume that $d(x)$ is bounded in $\Omega$.

We define recursively
\begin{eqnarray}
&& X_1(t) = (1- \log t)^{-1},\quad t\in (0,1],\nonumber \\
&& X_i(t) = X_{1}(X_{i-1}(t)), ~~~~~~~i=2,3,\ldots\, , t\in (0,1].
\label{va}
\end{eqnarray}
These are iterated logarithmic functions that vanish at an
increasingly slow rate at $t=0$ and satisfy $X_i(1)=1$. 
Given an integer $m \geq 1$ we define
\begin{equation}
 \eta_m(t) =  \sum_{i=1}^{m} X_1(t) \ldots X_i(t),  \qquad
 \zeta_m(t)= \sum_{i=1}^{m} X_1^2(t) \ldots X_i^2(t).
\label{oui}
 \end{equation}
We note that $\lim_{t\to 0}\eta_m(t)=\lim_{t\to 0}\zeta_m(t)=0$.
Now, it has been shown \cite{BFT2} that both series in (\ref{oui})
converge for any $t\in (0,1)$. This allows us to also introduce the functions $\eta_{\infty}$
and $\zeta_{\infty}$ as the infinite series.

We fix a parameter $s\in\R$ and we assume that the following inequality holds in the distributional sense:
\[p \neq k+s, \qquad
 (k+s-p)(d\Delta d- k+1)\geq 0
\quad\mbox{ in $\Omega  \setminus K$.} 
\]
For a detailed discussion of this condition we refer to
\cite{BFT1}. Here we simply note that it is satisfied in the
following two important cases: (i) it is satisfied as an equality
if $k=N$ and $K$ consists of single point and (ii) it is also
satisfied if $K=\partial\Omega$ (so $k=1$), $s+1-p<0$ and $\Omega$ is convex.

Our first theorem involves the functions
$\eta_{m}=\eta_{m}(d(x)/D)$ and $\zeta_{m}=\zeta_{m}(d(x)/D)$, $x\in\Omega$, for
a large enough parameter $D>0$. In any case $D$ will be large enough so that the quantity
$1+\alpha\eta_{m} +\beta\eta_{m}^2+\gamma\zeta_{m}$ is positive in $\Omega$.
We also set
\begin{equation}
H = \frac{k+s-p}{p}.
\label{eq:eta}
\end{equation}

\begin{theorem}
{\bf (weighted improved Hardy inequality)}
Let $p>1$ and $m\in\N\cup\{\infty\}$. Let $\Omega$ be a domain in $\R^N$ and $K$
a piecewise smooth surface of codimension $k$, $k=1,\ldots,N$.
Suppose that $p\neq k+s$, that $\sup_{x\in\Omega}d(x) < \infty$ and that 
\begin{equation}
 (k+s-p)(d\Delta d- k+1)\geq 0
\quad\mbox{ in $\Omega  \setminus K$.} 
\label{qqq}
\end{equation}
Also, let $\alpha,\beta,\gamma\in\R$ be fixed.
Then there exists a positive constant
$D_0\geq \sup_{x\in\Omega}d(x)$ such that for any $D \geq  D_0$
and all  $u \in C^{\infty}_c(\Omega\setminus K)$ there holds
\begin{eqnarray*}
&&\int_{\Omega}d^{s}(1+\alpha\eta_{m} +\beta\eta_{m}^2+\gamma\zeta_{m})|\nabla u|^pdx \geq
 |H|^p\int_{\Omega}d^{s-p}|u|^pdx +\\
&& \qquad\qquad +|H|^p\alpha\int_{\Omega}d^{s-p}\eta_{m} |u|^pdx
+\biggl( |H|^p\beta +\frac{|H|^{p-2}H\alpha}{2}\biggr)\int_{\Omega}d^{s-p}\eta_{m}^2|u|^pdx \\
&&\qquad\qquad +\biggl(\frac{p-1}{2p}|H|^{p-2} +\frac{|H|^{p-2}H\alpha}{2}
+ |H|^p\gamma\biggr)
\int_{\Omega}d^{s-p}\zeta_{m} |u|^pdx\, ,
\label{1.1}
\end{eqnarray*}
where $\eta_{m}=\eta_{m}(d(x)/D)= \sum_{i=1}^{m} X_1(d(x)/D) \ldots X_i(d(x)/D) $
and $\zeta_{m}=\zeta_{m}(d(x)/D)=\sum_{i=1}^{m} X_1^2(d(x)/D)\ldots X_i^2(d(x)/D)$.
\label{thma}
\end{theorem}
We note that the special case $s=\alpha=\beta=\gamma=0$ has been proved in \cite{BFT2}.

To state our next theorem we define the constant
\begin{equation}
Q=\frac{(p-1)k(k-2p)}{p^2}.
\label{qb}
\end{equation}

\begin{theorem}
{\bf (improved Rellich inequality I)}
Let $p>1$. Let $\Omega$ be a domain in $\R^N$ and $K$
a piecewise smooth surface of codimension $k$, $k=1,\ldots,N$.
Suppose that  $\sup_{x\in\Omega}d(x) < \infty$.
Suppose also that $k>2p$ and that
\[
d\Delta d-k+1\geq 0 \quad , \quad \mbox{in $\Omega\setminus K$}
\]
in the distributional sense. Then there exists a positive constant
$D_0\geq \sup_{x\in\Omega}d(x)$ such that for any $D \geq  D_0$
and all  $u \in C^{\infty}_c(\Omega\setminus K)$ there holds
\begin{eqnarray}
&&\int_{\Omega}|\Delta u|^pdx \geq
Q^p\int_{\Omega}\frac{|u|^p}{d^{2p}}dx +\label{1.2}\\
&& \qquad\qquad +\frac{p-1}{2p^3}|Q|^{p-2}\biggl\{k^2(p-1)^2+(k-2p)^2 \biggr\}\sum_{i=1}^{\infty}\int_{\Omega}
\frac{|u|^p}{d^{2p}}X_1^2X_2^2\ldots X_i^2 dx,
\nonumber
\end{eqnarray}
where $X_j=X_j(d(x)/D)$.
\label{thmb}
\end{theorem}

It is remarkable that the geometric assumption of this Theorem 2 only involves $\Delta d$, as in the case of Theorem 1, and not
higher-order derivatives of $d$ as one might expect.
The above theorem does not cover the important case $k=1$ which corresponds
to $d(x)=\dist(x,\partial\Omega)$. This is done in the following theorem for the case $p=2$.

\begin{theorem}
{\bf (improved Rellich inequality II)}
Let $\Omega$ be convex
and such that $d(x)\!\! :=$ $\dist(x,\partial\Omega)$ is bounded in $\Omega$.
Then there exists a positive constant
$D_0\geq \sup_{x\in\Omega}d(x)$ such that for any $D \geq  D_0$
and all  $u \in C^{\infty}_c(\Omega)$ there holds
\begin{eqnarray}
\ia && \int_{\Omega}(\Delta u)^2dx\geq \frac{1}{4}\int_{\Omega}\frac{|\nabla u|^2}{d^2}dx
+\frac{1}{4}\sum_{i=1}^{\infty}\int_{\Omega}
\frac{|\nabla u|^2}{d^2}X_1^2X_2^2\ldots X_i^2 dx \; , \\
\ib&& \int_{\Omega}(\Delta u)^2dx \geq
\frac{9}{16}\int_{\Omega}\frac{u^2}{d^4}dx +
\frac{5}{8}\sum_{i=1}^{\infty}\int_{\Omega}
\frac{u^2}{d^4}X_1^2X_2^2\ldots X_i^2 dx,
\label{sn}
\end{eqnarray}
where $X_j=X_j(d(x)/D)$.
\label{thmc}
\end{theorem}

In our last theorem we prove the optimality of the constants
appearing in Theorems 2 and 3 above. In a similar manner one can prove
the optimality of the constants in Theorem 1; we omit the proof
since it follows very closely the proof of \cite[Proposition
3.1]{BFT2}. Anyway, we note that in some particular cases the
optimality of Theorem 1 follows indirectly from the optimality of
Theorems 2 and 3, which we do prove. In relation to Theorem 4 see also the remark at the end of the paper.

We define
\[
J_0[u] =\int_{\Omega} |\Delta u|^p dx - |Q|^p
 \int_{\Omega} \frac{|u|^p}{d^{2p}} dx
\]
and for $m\in\N$,
\begin{eqnarray*}
J_m[u] &=&\int_{\Omega} |\Delta u|^p dx - |Q|^p
 \int_{\Omega} \frac{|u|^p}{d^{2p}} dx -  \\
&&- \frac{p-1}{2p^3}|Q|^{p-2}\biggl\{k^2(p-1)^2+(k-2p)^2 \biggr\}
 \sum_{i=1}^{m}\int_{\Omega} \frac{|u|^p}{d^{2p}} X_1^2 X_2^2\ldots X_i^2 dx.
\end{eqnarray*}
Our next theorem reads:
\begin{theorem}
Let $p>1$. Let $\Omega$ be a domain in $\R^N$. $\ia$ If $\, 2\leq k\leq N-1$ then we
take $K$ to be a piecewise smooth surface of codimension $k$ and assume $K\cap\Omega\neq
\emptyset$; $\ib$ if $k=N$ then we take $K=\{0\}\subset\Omega$; $\ic$ if $k=1$ then we assume
$K =\partial \Omega$. For any $D\geq \sup_{\Omega} d(x)$ we have
\begin{eqnarray*}
\ia && \inf_{C^{\infty}_c(\Omega\setminus K)}\frac{\int_{\Omega} |\Delta u|^p dx}{ \int_{\Omega}
\frac{|u|^p}{d^{2p}} dx} \leq |Q|^p; \\
\ib &&  \inf_{C^{\infty}_c(\Omega\setminus K)}\frac{J_{m-1}[u]}{\int_{\Omega} \frac{|u|^p}{d^{2p}} X_1^2 X_2^2
\ldots X_m^2}\leq  \frac{p-1}{2p^3}|Q|^{p-2}\biggl\{k^2(p-1)^2+(k-2p)^2 \biggr\},\; m\geq 1 .
\end{eqnarray*}
where $X_j=X_j(d(x)/D)$.
\label{thm:meli}
\end{theorem}

It follows in particular that all constants in Theorem \ref{thmb}
and Theorem \ref{thmc} (ii) are sharp. The sharpness of Theorem
\ref{thmc} (i) follows implicitly from the sharpness of \ref{thmc}
(ii).

\setcounter{equation}{0}
\section{Series expansion for weighted Hardy inequality}
In  this section we give the proof of Theorem \ref{thma}. We note that in the special
case $s=\alpha=\beta=\gamma=0$ the theorem has already been proved in \cite{BFT2}.
In the sequel we shall repeatedly use the differentiation rule
\begin{equation}
 \frac{d}{dt} X_i^{\beta}(t)=\frac{\beta}{t}X_1X_2\ldots X_{i-1}X_i^{1+\beta},
\qquad \beta\neq 0\, ,
\label{eq:3.52}
\end{equation}
which is easily proved by induction.

{\em Proof of Theorem \ref{thma}.}
We set for simplicity $\psi=(1+\alpha\eta_{m}+\beta\eta_{m}^2+\gamma\zeta_{m})$.
If $T$ is  a vector field in $\Omega$, then,
for any  $u \in C_c^{\infty}(\Omega \setminus K)$ we first
integrate by parts and then use  Young's inequality to obtain
\begin{eqnarray*}
\int_{\Omega} {\rm div}\, T ~  |u|^p dx
& \leq & p \int_{\Omega} | T| | \nabla u| |u|^{p-1} dx \\
& \leq &
  \int_{\Omega}  d^{s}\psi|\nabla u |^p dx  + (p-1)
 \int_{\Omega} d^{-\frac{s}{p-1}}|T|^{\frac{p}{p-1}}\psi^{-\frac{1}{p-1}}|u|^{p} dx\, ,
 \end{eqnarray*}
and thus conclude that
\begin{equation}
\int_{\Omega}d^{s} \psi | \nabla u |^p  dx  \geq
  \int_{\Omega} ( {\rm div}  T - (p-1)d^{-\frac{s}{p-1}}
|T|^{\frac{p}{p-1}}\psi^{-\frac{1}{p-1}}) |u|^{p} dx.
\label{2.21}
\end{equation}
We recall that $H=(k+s-p)/p$ and define
\[
T(x) = H  |H|^{p-2}d^{s+1-p}(x) \nabla d(x)\left( 1 +  (\alpha+\frac{p-1}{pH})  \eta_{m}(d(x)/D)
 + B \eta_{m}^2(d(x)/D) \right).
\]
where $D\geq\sup_{\Omega}d(x)$ and
$B\in\R$ is a free parameter to be chosen later. In any case, once $B$
is chosen, $D$ will be large enough so
that the quantity $1 +  (\alpha +\frac{p-1}{pH}) \eta_{m} (d/D) +B \eta_{m}^2(d/D)$ is positive
on $\Omega$. Note that $T$
is singular on $K$, but since  $u \in C_c^{\infty}(\Omega \setminus K)$
all previous calculations are legitimate.
In view of (\ref{2.21}), to prove the theorem it is enough to show that there exists $D_0\geq \sup_{\Omega}d(x)$
such that for $D\geq D_0$
\begin{eqnarray}
&&{\rm div}  T - (p-1)d^{-\frac{s}{p-1}}
|T|^{\frac{p}{p-1}}\psi^{-\frac{1}{p-1}}  -
d^{s-p}\Biggl\{ |H|^p
+|H|^p\alpha\eta_{m}  \label{eq:song} \\
&& \left.  +\biggl( |H|^p\beta +\frac{|H|^{p-2}H\alpha}{2}\biggr)\eta_{m}^2
+\biggl( \frac{p-1}{2p}|H|^{p-2} + \frac{|H|^{p-2}H\alpha}{2}+|H|^p\gamma\biggr) \zeta_{m}\right\}\geq 0
\nonumber
\end{eqnarray}
for all $x\in\Omega$.

To compute ${\rm div} T$ we shall need to differentiate $\eta_{m}(d/D)$. For this we
note that (\ref{eq:3.52}) easily implies
\[
\eta_{m}'(t) = \frac{1}{t} \left( X_1^2  + ( X_1^2 X_2 + X_1^2 X_2^2 )
 + \cdots +( X_1^2 X_2  \ldots X_m + \cdots + X_1^2  \ldots X_m^2)
\right),
\]
from which follows that
\begin{equation}
t \eta_{m}^{'} (t) = \frac{1}{2}  \zeta_{m}(t)
+ \frac{1}{2}  \eta_{m}^{2} (t).
\label{2.23}
\end{equation}
We also define $\theta_{m}$ on $(0,1)$ by
\[ \zeta_{m}'(t)=\frac{\theta_{m}(t)}{t},\]
and,
for simplicity, we set $A=\alpha+(p-1)/(pH)$ so that $T=|H|^{p-2}Hd^{s+1-p}\nabla d(1+A\eta_{m}+
B\eta_{m}^2)$.
We think of $\eta_{m}$ as an independent variable, which we may assume to be small by taking
$D$ large
enough. Simple computations together with assumption (\ref{qqq}) and the fact that $|\nabla d|=1$ give
\begin{eqnarray}
&&\diver T \\
&=& d^{s-p}\biggl\{p|H|^p+p|H|^pA\eta_{m}+
(p|H|^pB+|H|^{p-2}H\frac{A}{2})\eta_{m}^2 +|H|^{p-2}H\frac{A}{2}\zeta_{m} +\nonumber\\
&&\quad +|H|^{p-2}HB(\eta_{m}^3+\eta_{m}\zeta_{m}) \biggr\}\!+\!|H|^{p-2}Hd^{s-p}(d\Delta d+ s-p+1)
(1+A\eta_{m}+B\eta_{m}^2) \nonumber \\
&\geq&d^{s-p}\biggl\{p|H|^p+p|H|^pA\eta_{m}+
(p|H|^pB+|H|^{p-2}H\frac{A}{2})\eta_{m}^2 +|H|^{p-2}H\frac{A}{2}\zeta_{m} +\nonumber\\
&&\quad +|H|^{p-2}HB(\eta_{m}^3+\eta_{m}\zeta_{m}) \biggr\}+|H|^{p-2}Hd^{s-p}(d\Delta d-k+1)
(1+A\eta_{m}+B\eta_{m}^2) \nonumber \\
&\geq&d^{s-p}\biggl\{p|H|^p+p|H|^pA\eta_{m}+
(p|H|^pB+|H|^{p-2}H\frac{A}{2})\eta_{m}^2 +|H|^{p-2}H\frac{A}{2}\zeta_{m} +\nonumber\\
&&\quad +|H|^{p-2}HB(\eta_{m}^3+\eta_{m}\zeta_{m}) \biggr\}
\label{eq:ath}
\end{eqnarray}
Moreover, since $|\nabla d|=1$, Taylor's expansion gives
\begin{eqnarray}
|T|^{\frac{p}{p-1}}&=&|H|^pd^{\frac{(s+1-p)p}{p-1}}(1+A\eta_{m} +B\eta_{m}^2)^{\frac{p}{p-1}}\nonumber \\
&=&|H|^pd^{\frac{(s+1-p)p}{p-1}}\Biggl\{1+\frac{pA}{p-1}\eta_{m}+\biggl(\frac{pB}{p-1}+\frac{pA^2}{2(p-1)^2}\biggr)
\eta_{m}^2   \nonumber \\
&&+ \biggl(\frac{pAB}{(p-1)^2}-\frac{p(p-2)A^3}{6(p-1)^3}\biggr)\eta_{m}^3 +O(\eta_{m}^4)
\Biggr\}
\label{eq:ath1}
\end{eqnarray}
and also
\begin{eqnarray}
\psi^{-\frac{1}{p-1}}&=&1-\frac{\alpha}{p-1}\eta_{m} +\biggl(-\frac{\beta}{p-1}+
\frac{p\alpha^2}{2(p-1)^2}\biggr)\eta_{m}^2-\frac{\gamma}{p-1}\zeta_m\nonumber \\
&&+\biggl(\frac{p\alpha\beta}{(p-1)^2}-\frac{p(2p-1)\alpha^3}{6(p-1)^3}\biggr)\eta_{m}^3
+\frac{p\alpha\gamma}{(p-1)^2}\eta_{m}\zeta_{m}+O(\eta_{m}^4).
\label{eq:ath2}
\end{eqnarray}
Using (\ref{eq:ath}), (\ref{eq:ath1}) and (\ref{eq:ath2}) we see that the LHS of (\ref{eq:song}) is greater than or equal
to $d^{s-p}$ times a linear combination of powers of $\eta_{m}$, $\zeta_{m}$ and $\theta_{m}$
plus $O(\eta_{m}^4)$.
Recalling that $A=\alpha+(p-1)/(pH)$, we easily see that the constant term and the coefficients of $\eta_{m}$, $\eta_{m}^2$ and
$\zeta_{m}$ vanish, independently of the choice of the parameter $B$. The remaining
two coefficients, that is the coefficients of $\eta_{m}^3$ and $\eta_{m}\zeta_{m}$ are, respectively,
\[
\frac{(p-1)\alpha}{2pH^2}+\frac{\beta}{H}+\frac{(p-2)(p-1)}{6p^2H^3}\qquad  , \qquad \frac{B+\gamma}{H}.
\]
Since $\zeta_{m}\leq\eta_{m}^2\leq m\zeta_{m}$, we conclude that
taking $B$ to be large and positive (if $H>0$) or large and negative (if $H<0$), inequality
(\ref{eq:song}) is satisfied provided $\eta_{m}$ is small enough, which amounts to $D$ being large enough.
This completes the proof of the theorem. $\hfill //$

In the proof of Theorem \ref{thmb} we are going to use the last
theorem in the following special case which corresponds to taking $p=2$ and $s=-2q+2$ : \nl
{\bf Special case.} Assume that $k\neq 2q$ and that $(k-2q)(d\Delta d-k+1)\geq 0$ on $\Omega\setminus K$. Then for $D$ large enough
there holds
\begin{eqnarray}
&&\int_{\Omega}d^{-2q+2}(1+\alpha\eta_{m} +\beta\eta_{m}^2+\gamma\zeta_{m})|\nabla u|^2dx \geq
 \frac{(k-2q)^2}{4}\int_{\Omega}d^{-2q}u^2dx +\nonumber \\
&& \qquad +\frac{(k-2q)^2\alpha}{4}\int_{\Omega}d^{-2q}\eta_m u^2dx
+\biggl( \frac{(k-2q)^2\beta}{4} +\frac{(k-2q)\alpha}{4}\biggr)\int_{\Omega}d^{-2q}\eta_{m}^2u^2dx \nonumber \\
&&\qquad +\biggl(\frac{1}{4} +\frac{(k-2q)\alpha}{4}+\frac{(k-2q)^2\gamma}{4} \biggr)
 \int_{\Omega}d^{-2q}\zeta_{m}u^2dx
\label{eq:nd}
\end{eqnarray}
for all $u\in C^{\infty}_c(\Omega\setminus K)$.

\section{The improved Rellich inequality}

In this section we are going to prove Theorems \ref{thmb} and \ref{thmc} as well
as the corresponding optimality theorem. We begin with the following lemma where, we note,
$\Delta\phi$ is to be understood in the distributional sense.
\begin{lemma}
For any locally bounded function $\phi$ with $|\nabla u|\in L^2_{loc}(\Omega\setminus K)$ we have
\begin{equation}
\int_{\Omega} |\Delta u|^pdx \geq p(p-1)
\int_{\Omega}\phi |u|^{p-2}|\nabla u|^2dx -
\int_{\Omega}\Bigl( \Delta\phi +(p-1)|\phi|^{\frac{p}{p-1}}\Bigr)
|u|^pdx\, ,
\label{eq:vf}
\end{equation}
for all $u\in C^{\infty}_c(\Omega\setminus K)$.
\label{lem:vf}
\end{lemma}
{\em Proof.} Given $u\in C^{\infty}_c(\Omega\setminus K)$ we have
\begin{eqnarray*}
-\int_{\Omega}\Delta\phi |u|^pdx&=&p\int_{\Omega}\nabla\phi\cdot
(|u|^{p-2}u\nabla u)dx \\
&=&-p\int_{\Omega}\phi |u|^{p-2}u\Delta udx
-p(p-1)\int_{\Omega}\phi |u|^{p-2}|\nabla u|^2 dx \\
&\leq& p\left(\frac{p-1}{p}\int_{\Omega}|\phi|^{\frac{p}{p-1}}|u|^pdx
+\frac{1}{p}\int_{\Omega}|\Delta u|^pdx\right)- \\
&&\qquad -p(p-1)\int_{\Omega}\phi |u|^{p-2}|\nabla u|^2 dx.\\
\end{eqnarray*}
which is (\ref{eq:vf}).

{\em Proof of Theorem \ref{thmb}.} Let $m\in\N$ be fixed and let
$\eta_{m}$ and $\zeta_{m}$ be as in (\ref{oui}). We apply (\ref{eq:vf}) with
$\phi(x)=\lambda d(x)^{-2p+2}(1+\alpha\eta_{m}+\beta\eta_{m}^2)$, $\lambda>0$, where, as always,
$\eta_{m}=\eta_{m}(d(x)/D)$ and $D$ is yet to be determined. We thus obtain
\begin{equation}
\int_{\Omega}|\Delta u|^pdx \geq T_1+T_2+T_3
\label{eq:rad}
\end{equation}
where
\begin{eqnarray*}
&& T_1=p(p-1)\int_{\Omega}\phi |u|^{p-2}|\nabla u|^2dx\, , \\
&& T_2=- \int_{\Omega}\Delta\phi |u|^p dx\, ,\\
&& T_3 =-(p-1)\int_{\Omega}|\phi|^{\frac{p}{p-1}}|u|^pdx\, .
\end{eqnarray*}
To estimate $T_1$ we set $v=|u|^{p/2}$ and apply (\ref{eq:nd}) for $q=p$,
\begin{eqnarray}
T_1&=&\frac{4(p-1)\lambda}{p}\int_{\Omega}d^{-2p+2}(1+\alpha\eta_{m}+\beta\eta_{m}^2)
|\nabla v|^2dx \nonumber\\
&\geq& \frac{4(p-1)\lambda}{p}\int_{\Omega}d^{-2p}
\left\{\frac{(k-2p)^2}{4}
+\frac{(k-2p)^2\alpha}{4}\eta_{m} +\right. \nonumber \\
&&\left.\biggl(\frac{1}{4}+
\frac{(k-2p)\alpha}{4}\biggr)\zeta_{m} +
\biggl(\frac{(k-2p)\alpha}{4}+\frac{(k-2p)^2\beta}{4}\biggr)
\eta_{m}^2 \right\}|u|^pdx
\label{t1}
\end{eqnarray}
To estimate $T_2$ we first note that
\[\nabla\phi=\lambda d^{-2p+1}\left\{-2(p-1)(1+\alpha\eta_{m}+\beta\eta_{m}^2)
+\frac{\alpha}{2}(\eta_{m}^2+\zeta_{m}^2)+\beta(\eta_{m}^3+\eta_{m}\zeta_{m})\right\}\nabla d\]
and hence compute
\begin{eqnarray*}
-\Delta\phi &=&\lambda d^{-2p}(-2p+1+d\Delta d)\biggl\{2(p-1)(1+\alpha\eta_{m}+
\beta\eta_{m}^2)-\frac{\alpha}{2}(\eta_{m}^2+\zeta_{m})-\\
&&\hspace{8cm}-\beta(\eta_{m}^3+\eta_{m}\zeta_{m})\biggr\}+\\
&&+\lambda d^{-2p}\biggl\{(p-1)\alpha(\eta_{m}^2+\zeta_{m})-\beta(\eta_{m}^3+\eta_{m}\zeta_{m})- \\
&&\hspace{6cm}-\frac{\alpha}{2}(\eta_{m}^3+\eta_{m}\zeta_{m}+\theta_{m})+O(\eta_{m}^4)\biggr\}
\end{eqnarray*}
Using the geometric assumption $d\Delta d-k+1\geq 0$
and collecting similar terms we conclude that
\begin{eqnarray}
T_2&\geq&\lambda \int_{\Omega}d^{-2p}\Biggl\{ 2(p-1)(k-2p)+
2(p-1)(k-2p)\alpha\eta_{m} +
\nonumber\\
&&+\biggl(2(p-1)(k-2p)\beta+\frac{-k+4p-2}{2}\alpha\biggr)\eta_{m}^2 \label{t2}\\
&&+\frac{-k+4p-2}{2}\alpha\zeta_{m} -\biggl( (k-2p+1)\beta+\frac{\alpha}{2}\biggr)
(\eta_{m}^3+\eta_{m}\zeta_{m})-\frac{\alpha}{2}\theta_{m} +O(\eta_{m}^4)\Biggr\}dx.
\nonumber
\end{eqnarray}
From Taylor's theorem we have
\begin{eqnarray*}
(1+\alpha\eta_{m}+\beta\eta_{m}^2)^{\frac{p}{p-1}}&=&1+\frac{p\alpha}{p-1}\eta_{m} +
\biggl(\frac{p\beta}{p-1}+\frac{p\alpha^2}{2(p-1)^2}\biggr)\eta_{m}^2 +\\
&& +\biggl(\frac{p\alpha\beta}{(p-1)^2}+\frac{p(2-p)\alpha^3}{6(p-1)^3}
\biggr)\eta_{m}^3 +O(\eta_{m}^4)
\end{eqnarray*}
from which follows that
\begin{eqnarray}
T_3&=&-(p-1)|\lambda|^{\frac{p}{p-1}}\int_{\Omega}d^{-2p}
\Biggl\{1+\frac{p\alpha}{p-1}\eta_{m} +
\biggl(\frac{p\beta}{p-1}+\frac{p\alpha^2}{2(p-1)^2}\biggr)\eta_{m}^2 +\nonumber\\
&&\hspace{3cm}+\biggl(\frac{p\alpha\beta}{(p-1)^2}+\frac{p(2-p)\alpha^3}{6(p-1)^3}
\biggr)\eta_{m}^3 +O(\eta_{m}^4)\Biggr\}dx.
\label{t3}
\end{eqnarray}
Using the above estimates on $T_1,T_2$ and $T_3$ and going back to (\ref{eq:rad})
we obtain the inequality
\begin{equation}
\int_{\Omega}|\Delta u|^pdx\geq \int_{\Omega}d^{-2p}V|u|^pdx
\label{eq:v}
\end{equation}
where the potential $V$ has the form
\[ V(x)=r_0+r_1\eta_{m} +r_2\eta_{m}^2 +r_2'\zeta_{m} +r_3\eta_{m}^3+r_3'\eta_{m}\zeta_{m}+r_3''\theta_{m}.\]
We compute the coefficients $r_i,r_i'$ by adding the corresponding coefficients from
(\ref{t1}), (\ref{t2}) and (\ref{t3}). We ignore for now the coefficients
of the third-order terms. For the others we find
\begin{eqnarray*}
&&r_0=(p-1)\biggl(\frac{k(k-2p)}{p}\lambda-|\lambda|^{\frac{p}{p-1}}\biggr)\\
&&r_1=\frac{(p-1)k(k-2p)}{p}\alpha\lambda -p\alpha|\lambda|^{\frac{p}{p-1}}\\
&&r_2= \frac{pk-2k+2p}{2p}\alpha\lambda
+\frac{(p-1)k(k-2p)}{p}\beta\lambda -(p-1)\biggl(\frac{p\beta}{p-1}+
\frac{p\alpha^2}{2(p-1)^2}\biggr)|\lambda|^{\frac{p}{p-1}}\\
&&r_2'=\biggl(\frac{p-1}{p}+\frac{pk-2k+2p}{2p}\alpha\biggr)\lambda
\end{eqnarray*}
We now make a specific choice for $\alpha$ and $\lambda$. We recall that
$Q=(p-1)k(k-2p)/p^2$,
and choose
\[ \lambda=Q^{p-1}\qquad , \qquad \alpha =\frac{(p-1)(pk-2k+2p)}{p^2Q} .\]
We then have $r_0=Q^p$, $r_1=r_2=0$, irrespective of the value of $\beta$.
We also have
\[ r_2'=\frac{p-1}{p}Q^{p-2}\biggl(Q+\frac{(pk-2k+2p)^2}{2p^2}\biggr).\]
Substituting these values in (\ref{eq:v}) we thus obtain
\begin{eqnarray*}
\int_{\Omega}|\Delta u|^pdx \!\!\!&\geq&\! Q^p\int_{\Omega}d^{-2p}|u|^pdx
+\frac{p-1}{p}Q^{p-2}\biggl(Q+\frac{(pk-2k+2p)^2}{2p^2}\biggr)
\int_{\Omega}d^{-2p}\zeta_{m} |u|^pdx \\
&& +\int_{\Omega}d^{-2p}(r_3\eta_{m}^3+r_3'\eta_{m}\zeta_{m}+r_3''\theta_{m}
+O(\eta_{m}^4))|u|^pdx
\end{eqnarray*}
We still have not imposed any restriction on $\beta$. We now
observe that $r_3'$ and $r_3''$ are independent of $\beta$, while
$r_3=c_1\beta+c_2$ with $c_1=Q^{p-1}((2k)/p-2k+2p-3)<0$. Hence,
since the functions $\eta_{m}^3$, $\eta_{m}\zeta_{m}$ and
$\theta_{m}$ are comparable in size to each other, the integral is
made positive by choosing $\beta$ to be large and negative and
$\eta_{m}$ small enough, which amounts to $D$ being large enough.
Hence we have proved that for $D\geq D_0$ there holds
\begin{eqnarray*}
&&\int_{\Omega}|\Delta u|^pdx \geq
Q^p\int_{\Omega}\frac{|u|^p}{d^{2p}}dx +\\
&& \qquad\qquad +\biggl(\frac{p-1}{p}Q^{p-1}+\frac{p-1}{2p}Q^{p-2}R^2\biggr)\sum_{i=1}^{m}\int_{\Omega}
\frac{|u|^p}{d^{2p}}X_1^2X_2^2\ldots X_i^2 dx,
\end{eqnarray*}
where $X_j=X_j(d(x)/D)$. This concludes the proof of the theorem. $\hfill //$

{\bf Remark.}  Let us mention here that in the proofs of Theorems
\ref{thma} and \ref{thmb} we did not use at any point the
assumption that $k$ is the codimension of the set $K$. Indeed, a
careful look at the two proofs shows that $K$ can be any closed
set such that $\dist(x,K)$ is bounded in $\Omega$ and for which
the condition  $d\Delta d -k+1\geq 0$ or $\leq 0$ is satisfied; the proof does not even require $k$ to
be an integer. Of course, the natural realizations of these
conditions are that $K$ is smooth and $k=\codim(K)$. However, the
argument also applies in the case where $K$ is a union of sets of
different codimensions; see \cite{BFT1}.

We next prove Theorem \ref{thmc}.

{\em Proof of Theorem \ref{thmc}.} We note that the convexity of $\Omega$ implies
that $\Delta d\leq 0$ on $\Omega$ in the distributional sense \cite[Theorem 6.3.2]{EG}.
Now, let $u\in C^{\infty}_c(\Omega)$ be given. Applying Theorem \ref{thma}
(with $k=1$, $p=2$, $s=0$ and $\alpha=\beta=\gamma=0$) to the
partial derivatives $u_{x_i}$ we have
\begin{eqnarray}
\int_{\Omega}(\Delta u)^2dx &=&\sum_{i=1}^n\int_{\Omega}|\nabla u_{x_i}|^2dx\nonumber\\
&\geq&\sum_{i=1}^n\Biggl\{\frac{1}{4}\int_{\Omega}\frac{u_{x_i}^2}{d^2}dx
+\frac{1}{4}\int_{\Omega}\frac{u_{x_i}^2}{d^2}\zeta_{m} dx \nonumber\\
&=&\frac{1}{4}\int_{\Omega}\frac{|\nabla u|^2}{d^2}(1+\zeta_{m})dx
\Biggr\},
\label{ci}
\end{eqnarray}
for $D$ large enough, where $\zeta_{m}=\zeta_{m}(d(x)/D)$.
Applying Theorem \ref{thma} once more (this time with $k=1$, $p=2$, $s=-2$, $\alpha=\beta=0$
and $\gamma=1$)
we obtain
\begin{equation}
\int_{\Omega}\frac{|\nabla u|^2}{d^2}(1+\zeta_{m})dx \geq\frac{9}{4}\int_{\Omega}\frac{u^2}{d^4}dx
+\frac{5}{2}\int_{\Omega}\frac{u^2}{d^4}\zeta_{m} dx\, .
\label{ele}
\end{equation}
Combining (\ref{ci}) and (\ref{ele}) we obtain
\[
\int_{\Omega}(\Delta u)^2dx \geq
\frac{9}{16}\int_{\Omega}\frac{u^2}{d^4}dx
+\frac{5}{8}\int_{\Omega}\frac{u^2}{d^4}\zeta_{m} dx\, ,
\]
which is the stated inequality. $\hfill //$


We next give the proof of Theorem \ref{thm:meli}.
We recall that $\Omega$ is a domain in $\R^N$ and
that $K$ is a piecewise smooth surface of codimension $k$ such that $K\cap\Omega\neq
\emptyset$, unless $k=1$ in which case $K =\partial \Omega$.
All the calculations below are local, in a small ball of radius $\delta$, and indeed,
it would be enough to assume that $K$ has a smooth part.
We also note that for $k=N$ (distance from a point) the subsequent calculations are substantially
simplified, whereas for $k=1$ (distance from the boundary) one should
replace $B_{\delta}$ by $B_{\delta} \cap \Omega$. This last  change
entails some minor modifications, the arguments otherwise being  the same.

{\em Proof.} We shall only give the proof of $\ib$ since the proof of $\ia$ is much simpler.
For the proof we shall use some of the ideas and tools developed in \cite{BFT2}.
All our analysis will be local, say, in a fixed ball of $B(x_0,\delta)$
where $x_0\in K$ and $\delta$ is small, but fixed throughout th proof.
We therefore fix a smooth, non-negative  function $\phi$ such that $\phi(x)=1$ on $\{|x-x_0|<\delta/2\}$ and
$\phi(x)=1$ on $\{|x-x_0|>\delta\}$.
For given $\epsilon_0,\epsilon_1\ldots,\epsilon_m>0$ we then define the function
\begin{eqnarray*}
u&=&\phi d^{\frac{-k+2p+\epsilon_0}{p}}X_1^{\frac{-1+\epsilon_1}{p}}X_2^{\frac{-1+\epsilon_2}{p}}
\ldots X_m^{\frac{-1+\epsilon_m}{p}} \\
&=:&\phi v.
\end{eqnarray*}
A standard argument using cut-off functions shows that $u$ belongs in $W^{2,p}_0(\Omega\setminus K)$
and therefore is a legitimate test-function for the infimum above.
We intend to see how $J_{m-1}[u]$ behaves as the $\epsilon_i$'s tend to zero.
We shall not be interested in terms that remain bounded for small values of the
$\epsilon_i$'s. To distinguish such terms we shall need the following fact, cf. \cite[(3.8)]{BFT2}:
we have
\begin{eqnarray}
&&\int_{\Omega}\phi^p d^{-k+\beta_0}X_1^{1+\beta_1}(d/D)\ldots X_m^{1+\beta_m}(d/D)dx
<\infty \Longleftrightarrow \nonumber \\
&&\Longleftrightarrow\; \left\{
\begin{array}{ll}
& \beta_0>0 \\
\mbox{ or}&\mbox{$\beta_0=0$ and $\beta_1>0$}\\
\mbox{ or}&\mbox{$\beta_0=\beta_1=0$ and $\beta_2>0$}\\
& \cdots \\
\mbox{ or}&\mbox{$\beta_0=\beta_1=\ldots =\beta_{m-1}=0$ and $\beta_m>0$.}
\end{array}\right.
\label{eq:3.51}
\end{eqnarray}
Now, we have $\Delta u=\phi\Delta v+2\nabla\phi\cdot\nabla v+v\Delta\phi$ and hence, using the inequality
\begin{equation}
|a+b|^p\leq |a|^p +c(|a|^{p-1}|b|+|b|^p),
\label{eq:ab}
\end{equation}
we have
\begin{eqnarray*}
&& \int_{\Omega}|\Delta u|^pdx \\
&\leq& \int_{\Omega}\phi^p|\Delta v|^pdx +c \int_{\Omega}\biggl\{(\phi|\Delta v|)^{p-1}
(|\nabla\phi||\nabla v|+|v||\Delta\phi|)+(|\nabla\phi||\nabla v|+|v||\Delta\phi|)^p\biggr\}\\
&\leq& \int_{\Omega}\phi^p|\Delta v|^pdx +\int_{\Omega}\biggl( |\Delta v|^{p-1}(|\nabla v|+|v|)+(|\nabla v|+|v|)^p\biggr).
\end{eqnarray*}
The first integral involves $d$ to the power $-k+\epsilon_0/p$ (see below) and is therefore important.
On the other hand, all terms in the second integral involve $d$ to the power that is larger than $-k$
and in fact bounded away from $-k$, independently of $\epsilon_0$; hence
\begin{equation}
\int_{\Omega}|\Delta u|^pdx =\int_{\Omega}\phi^p |\Delta v|^p dx +O(1),
\label{eq:ps}
\end{equation}
where the $O(1)$ is uniform in all the $\epsilon_i$'s.

We next define the function
\[ g(t)=\epsilon_0 +(-1+\epsilon_1)X_1+(-1+\epsilon_2)X_1X_2+\cdots (-1+\epsilon_m)X_1X_2\ldots X_m,\quad
t>0, \]
where $X_i=X_i(t/D)$. We shall always think of $g(d(x))$ as a small quantity.
Recalling (\ref{eq:3.52}) one easily sees that
\begin{eqnarray}
\frac{d\eta_{m}}{dt}&=&\frac{1}{t}\sum_{1\leq i\leq j\leq m}(-1+\epsilon_j)X_1^2\ldots X_i^2X_{i+1}\ldots X_j
\nonumber \\
&=:&\frac{h(t)}{t}.
\label{eq:vas}
\end{eqnarray}
Also, for any $\beta$ there holds
\begin{equation}
\frac{d}{dt} \bigl(t^{\frac{\beta+\epsilon_0}{p}}X_1^{\frac{-1+\epsilon_1}{p}}X_2^{\frac{-1+\epsilon_2}{p}}\!\!
\ldots X_m^{\frac{-1+\epsilon_m}{p}} \bigr)=
t^{\frac{\beta-p+\epsilon_0}{p}}X_1^{\frac{-1+\epsilon_1}{p}}X_2^{\frac{-1+\epsilon_2}{p}}\!\!
\ldots X_m^{\frac{-1+\epsilon_m}{p}}[\frac{\beta}{p}+\frac{g(t)}{p}].
\label{eq:g}
\end{equation}
Applying (\ref{eq:g}) first for $\beta=k-2p$, then for $\beta=k-p$ and using (\ref{eq:vas}) we obtain
\[ \Delta v=d^{\frac{-k+\epsilon_0}{p}}X_1^{\frac{-1+\epsilon_1}{p}}X_2^{\frac{-1+\epsilon_2}{p}}
\ldots X_m^{\frac{-1+\epsilon_m}{p}}\Biggl\{\biggl( \frac{k-p}{p}-d\Delta d-\frac{g}{p}\biggr)
\biggl( \frac{k-2p}{p}-\frac{g}{p} \biggr) +\frac{h}{p}\Biggr\},
\]
where, here and below, we use $g$, $h$ and $X_i$ to denote $g(d(x))$, $h(d(x))$ and $X_i(d(x)/D)$.
Now, by \cite[Theorem 3.2]{AS}  we have $d\Delta d=k-1+O(d)$ as $d(x)\to 0$.
Hence, the expression in the braces equals
\[Q+\frac{R}{p}g-\frac{1}{p^2}g^2-\frac{1}{p}h +O(d) \mbox{ as } d(x)\to 0,\]
where $R=(2k-pk-2p)/p$.
The $O(d)$ gives a bounded contribution by an application of (\ref{eq:ab}) -- as was done earlier.
Hence (\ref{eq:ps}) gives
\begin{equation}
\int_{\Omega}|\Delta u|^pdx=\int_{\Omega}\phi^pd^{-k+\epsilon_0}X_1^{-1+\epsilon_1}\!\! \ldots X_m^{-1+\epsilon_m}
\biggl| Q+\frac{R}{p}g-\frac{1}{p^2}g^2-\frac{1}{p}h \biggr|^pdx +O(1).
\label{eq:ps1}
\end{equation}
To estimate this we take the Taylor's expansion of $|Q+t|^p$ about $t=0$. We obtain
\begin{eqnarray}
\int_{\Omega}|\Delta u|^pdx&=&\int_{\Omega}\phi^pd^{-k+\epsilon_0}X_1^{-1+\epsilon_1}\!\! \ldots X_m^{-1+\epsilon_m}
\Biggl\{ |Q|^p+|Q|^{p-2}QRg +  \nonumber \\
&&\quad +\biggl(-\frac{1}{p}|Q|^{p-2}Q+\frac{p-1}{2p}|Q|^{p-2}R^2\biggr)g^2- \nonumber \\
&&\quad - |Q|^{p-2}Q \zeta_{m} +O(g^3)+O(gh)+O(h^2)\Biggr\}dx +O(1).
\label{eq:ps2}
\end{eqnarray}
Using (\ref{eq:ab}) once again, it is not difficult to see that the terms $O(g^3)$, $O(gh)$ and $O(h^2)$ give a contribution
that is bounded uniformly in the $\epsilon_i$'s and can therefore be dropped.

At this point, and in order to simplify the notation, we introduce some auxiliary quantities.
For $0\leq i\leq j\leq m$ we define
\begin{eqnarray*}
&& A_0 =\int_{\Omega}\phi^p d^{-k+\epsilon_0}X_1^{-1+\epsilon_1}\!\!\!\!
\ldots X_m^{-1+\epsilon_m} dx \\
&& A_i =\int_{\Omega}\phi^p d^{-k+\epsilon_0}X_1^{1+\epsilon_1}\!\!\!\!\ldots
X_i^{1+\epsilon_i}X_{i+1}^{-1+\epsilon_{i+1}}\!\!\!\!\ldots X_m^{-1+\epsilon_m} dx \\
&& \Gamma_{0j}=\int_{\Omega}\phi^p d^{-k+\epsilon_0}X_1^{\epsilon_1}\!\!\!\!\ldots
X_i^{\epsilon_i}X_{i+1}^{-1+\epsilon_{i+1}}\!\!\!\!\ldots X_m^{-1+\epsilon_m} dx \\
&& \Gamma_{ij}=\int_{\Omega}\phi^p d^{-k+\epsilon_0}X_1^{1+\epsilon_1}\!\!\!\!\ldots
X_i^{1+\epsilon_i}X_{i+1}^{\epsilon_{i+1}}\!\!\!\!\ldots X_j^{\epsilon_j}
X_{j+1}^{-1+\epsilon_{j+1}}\!\!\!\!\ldots X_m^{-1+\epsilon_m} dx,
\end{eqnarray*}
with the convention that $\Gamma_{ii}=A_i$. It is then easily seen that
\begin{eqnarray*}
&&\int_{\Omega}\phi^p d^{-k+\epsilon_0}X_1^{-1+\epsilon_1}\!\!\!\!
\ldots X_m^{-1+\epsilon_m} g\,dx =\epsilon_0 A_0-\sum_{i=1}^m(1-\epsilon_i)\Gamma_{0i} \\
&&\int_{\Omega}\phi^p d^{-k+\epsilon_0}X_1^{-1+\epsilon_1}\!\!\!\!
\ldots X_m^{-1+\epsilon_m} g^2dx =\epsilon_0^2A_0+\sum_{i=1}^m(1-\epsilon_i)^2A_i
-2\epsilon_0\sum_{i=1}^m(1-\epsilon_i)\Gamma_{0i}+\\
&&\hspace{6.5cm}+2\sum_{i<j}(1-\epsilon_i)(1-\epsilon_j)\Gamma_{ij} \\
&&\int_{\Omega}\phi^p d^{-k+\epsilon_0}X_1^{-1+\epsilon_1}\!\!\!\!
\ldots X_m^{-1+\epsilon_m} h\, dx=-\sum_{i=1}^m(1-\epsilon_i)A_i -\sum_{i<j}(1-\epsilon_j)\Gamma_{ij}.
\end{eqnarray*}
(Here and below $\sum_{i<j}$ means $\sum_{1\leq i<j\leq m}$.)
Let us also define the constant
\[ P=-\frac{1}{p}|Q|^{p-2}Q+\frac{p-1}{2p}|Q|^{p-2}R^2 .\]
Going back to (\ref{eq:ps2}) and noting that
\begin{eqnarray*}
\int_{\Omega}\frac{|u|^p}{d^{2p}}dx &=&A_0  \\
\sum_{i=1}^{m-1}\int_{\Omega}\frac{|u|^p}{d^{2p}}X_1^2\ldots X_i^2dx&=&\sum_{i=1}^{m-1}A_i.
\end{eqnarray*}
we obtain
\begin{eqnarray}
J_{m-1}[u]&=&\! |Q|^{p-2}QR\biggl(\epsilon_0A_0-\sum_{i=1}^m(1-\epsilon_i)\Gamma_{0i}\biggr)+\nonumber \\
&&\!\! +P\biggl(\epsilon_0^2A_0+ \sum_{i=1}^m(1-\epsilon_i)^2A_i
-2\epsilon_0\sum_{i=1}^m(1-\epsilon_i)\Gamma_{0i}+2\sum_{i<j}(1-\epsilon_i)(1-\epsilon_j)\Gamma_{ij} \biggr)\nonumber \\
&&\! +|Q|^{p-2}Q\biggl(\sum_{i=1}^m(1-\epsilon_i)A_i +\sum_{i<j}(1-\epsilon_j)\Gamma_{ij}\biggr)
-G\sum_{i=1}^{m-1}A_i+O(1).
\label{eq:zeus}
\end{eqnarray}
Now, by \cite[p184]{BFT2},
\[\epsilon_0^2-2\epsilon_0\sum_{i=1}^m(1-\epsilon_i)\Gamma_{0i}=
\sum_{i=1}^m(\epsilon_i-\epsilon_i^2)A_i+\sum_{i<j}(2\epsilon_i-1)(1-\epsilon_j)
\Gamma_{ij} +O(1).\]
For the sake of simplicity, we set
\[
G=\frac{p-1}{2p^3}|Q|^{p-2}\biggl\{k^2(p-1)^2+(k-2p)^2 \biggr\}.
\]
One then easily sees that $P+|Q|^{p-2}Q=G$. Hence, collecting similar terms,
\begin{eqnarray*}
J_{m-1}[u]&=&|Q|^{p-2}QR(\epsilon_0A_0-\sum_{j=1}^m(1-\epsilon_j)\Gamma_{0j}
-G\biggl(\sum_{i=1}^m\epsilon_iA_i-\sum_{i<j}(1-\epsilon_j)\Gamma_{ij}
\biggr)+\\
&&+GA_m +O(1),
\end{eqnarray*}
where the $O(1)$ is uniform for small $\epsilon_i$'s.

Up to this point the parameters $\epsilon_0,\epsilon_1,\ldots,\epsilon_m$
where positive. We intend to take limits as they tend to zero in that order.
Due to (\ref{eq:3.51}), as $\epsilon_0\to 0$ all terms have
finite limits except those involving $A_0$ and $\Gamma_{0j}$ which, when viewed separately,
diverge. However a simple argument involving an integration by parts
(see \cite[(3.9)]{BFT2} shows that
\begin{equation}
\epsilon_0A_0-\sum_{j=1}^m(1-\epsilon_j)\Gamma_{0j}=O(1)
\label{eq:air}
\end{equation}
uniformly in $\epsilon_0,\ldots,\epsilon_m$. Hence, letting $\epsilon_0\to 0$
we conclude that
\[
J_{m-1}[u]=-G\biggl(\sum_{i=1}^m\epsilon_iA_i-\sum_{i<j}(1-\epsilon_j)\Gamma_{ij}
\biggr)+GA_m +O(1) \qquad (\epsilon_0=0)
\]
Now -- as was the case with (\ref{eq:air}) -- an integration by parts shows that (see \cite[(3.9)]{BFT2})
if $\epsilon_0 = \epsilon_1 =\ldots=\epsilon_{i-1}=0$, then
\begin{equation}
\epsilon_iA_i-\sum_{j=1}^m(1-\epsilon_j)\Gamma_{ij}=O(1).
\label{eq:air1}
\end{equation}
We now let $\epsilon_1\to 0$. Again, all terms have finite limits
except those involving $A_1$ and $\Gamma_{1j}$ which diverge.
Using (\ref{eq:air1}) we see that when combined these terms stay
bounded in the limit $\alpha_1\to 0$. We proceed in this way and
after letting $\epsilon_{m-1}\to 0$ we are left with
\[ J_{m-1}[u]=G(1-\epsilon_m)A_m +O(1)\qquad (\epsilon_0=\ldots=\epsilon_{m-1}=0).\]
Let us denote by $G'$ the infimum in the left-hand side of part (ii) of Theorem 4. We have thus proved that
\[ G'\leq \frac{G(1-\epsilon_m)A_m+O(1)}{A_m}  \qquad (\epsilon_0=\ldots=\epsilon_{m-1}=0).\]
Letting now $\epsilon_m\to 0$ we have $A_m\to +\infty$ (by \ref{eq:3.51})), and thus conclude that
$G'\leq G$, as required. $\hfill //$

{\bf Remark.} Slightly modifying the above argument one can also prove the optimality
of the power $X_m^2$ of the improved Rellich inequalities (\ref{1.2})
and (\ref{sn}). Namely, for any $\epsilon>0$ there holds
\[
\inf_{u \in C^{\infty}_c(\Omega\setminus K)}
 \frac{J_{m-1}[u]}
{\int_{\Omega} \frac{|u|^p}{d^{2p}}X_1^2 X_2^2
 \ldots X_m^{2-\epsilon}dx}=0\, .
\]

{\bf Acknowledgment} We acknowledge partial  support by the  RTN  European network
 Fronts--Singularities, HPRN-CT-2002-00274.



\end{document}